\documentclass[12pt]{article}

\usepackage{amsmath,amsthm,amssymb}

\usepackage{txfonts}	

\usepackage{color}
\definecolor{myColor}{gray}{0.975}

\newcommand{\create}[1]{\begin{#1}}
\newcommand{\delete}[1]{\end{#1}}

\newcommand{\Frac}[2]{\leavevmode\kern.10em\raise.50ex\hbox{\the\scriptfont0 $#1$}\kern+.20em
                         {\big /}\kern.20em\lower.80ex\hbox{\the\scriptfont0 $#2$}}

\newcommand{\Sum}[2]
{
	\lower.40ex\hbox{$
	\underset{#1}{\overset{#2}{\mathlarger{\mathlarger{\mathlarger \Sigma}}}}
	\hspace{.04in}$}
}
\newcommand{\cupEx}[2]
{
	\lower.40ex\hbox{$
	\underset{#1}{\overset{#2}{\mathlarger{\mathlarger{\mathlarger{ \cup}}}}}
	\hspace{.02in}$}
}
\newcommand{\capEx}[2]
{
	\lower.40ex\hbox{$
	\underset{#1}{\overset{#2}{\mathlarger{\mathlarger{\mathlarger{ \cap}}}}}
	\hspace{.02in}$}
}
\newcommand{\Prod}[2]
{
	\lower.40ex\hbox{$
	\underset{#1}{\overset{#2}{\mathlarger{\mathlarger{\mathlarger \Pi}}}}
	\hspace{.04in}$}
}
\newcommand{\Directsum}{\rotatebox[origin=c]{180}{\reflectbox{$\Pi$}}}
\newcommand{\Coprod}[2]
{
	\lower.40ex\hbox{$
	\underset{#1}{\overset{#2}{\mathlarger{\mathlarger{\mathlarger \Directsum}}}}
	\hspace{.04in}$}
}
\newcommand{\OPlus}[2]
{
	\lower.40ex\hbox{$
	\underset{#1}{\overset{#2}{\mathlarger{\mathlarger{\mathlarger \oplus}}}}
	\hspace{.02in}$}
}

\newcommand{\simto} {
	\rightarrow\hspace{-12.5pt}\raise4.pt\hbox{$\mathsmaller \sim$} \hspace{7.4pt}}
\newcommand{\simlongto} {\longrightarrow \hspace{-16.5pt}\raise4pt\hbox{$\sim$} \hspace{5.51pt}}
\newcommand{\simlongtoL}{
	\longrightarrow \hspace{-16.5pt}\raise4pt\hbox{$\sim$} \hspace{8.85pt}}

\newcommand{\coloniff}{\hspace{.1in}\raise0.06ex\hbox{:}\mspace{-10mu}\iff}

\newcommand{\sen}
{
	\noindent
	\hspace{-2em}
	\leavevmode
	\leaders\hrule depth-2.1pt height 3.3pt\hfill\kern0pt
	\hspace{-2em}
	\vspace{.5ex}
}

\newcommand{\senS}
{
	\noindent
	\hspace{-.6em}
	\leavevmode
	\leaders\hrule depth-2.1pt height 2.5pt\hfill\kern0pt
	\hspace{-.6em}
	\vspace{.5ex}
}

\newcommand{\Ker}{\text{Ker}\,}

\newcommand{\hs}[1]{\hspace{#1}}
\newcommand{\vs}[1]{\vspace{#1}}

\newcommand{\MARGIN}{\vspace{-1ex}}
\newcommand{\beginhspace}{\hspace{.03in}}

\newcommand{\TxtBgn}{\vs{.5ex}}
\newcommand{\TxtEnd}{\vs{.5ex}}

\newcommand{\NtnBgn}{\MARGIN\begin{Ntn} \normalfont \beginhspace}
\newcommand{\NtnEnd}{\hfill\end{Ntn}\MARGIN}

\newcommand{\DefBgn}{\MARGIN\begin{Def} \normalfont \beginhspace}
\newcommand{\DefEnd}{\hfill \end{Def}\MARGIN}

\newcommand{\PrpBgn}{\MARGIN\begin{Prop} \normalfont \beginhspace}
\newcommand{\PrpEnd}{\hfill \end{Prop}\MARGIN}

\newcommand{\ThmBgn}{\MARGIN\begin{Thm} \normalfont \beginhspace}
\newcommand{\ThmEnd}{\hfill \end{Thm}\MARGIN}

\newcommand{\CorBgn}{\MARGIN\begin{Cor} \normalfont \beginhspace}
\newcommand{\CorEnd}{\hfill \end{Cor}\MARGIN}

\newcommand{\LemBgn}{\MARGIN\begin{Lem} \normalfont \beginhspace}
\newcommand{\LemEnd}{\hfill \end{Lem}\MARGIN}

\newcommand{\ClmBgn}{\MARGIN\begin{Clm} \normalfont \beginhspace}
\newcommand{\ClmEnd}{\hfill \end{Clm}\MARGIN}

\newcommand{\PrfBgn}{\MARGIN\begin{Proof} \normalfont \beginhspace}
\newcommand{\PrfEnd}{\hfill\end{Proof}\MARGIN}

\newcommand{\RmkBgn}{\MARGIN\begin{Remark} \normalfont \beginhspace}
\newcommand{\RmkEnd}{\hfill \end{Remark}\MARGIN}

\newcommand{\ExpBgn}{\MARGIN\begin{Ex} \normalfont \beginhspace}
\newcommand{\ExpEnd}{\hfil \end{Ex}\MARGIN}

\newtheorem{Def}{Definition.}[section]
\newtheorem{Ntn}[Def]{Notation.}
\newtheorem{Thm}[Def]{Theorem.}
\newtheorem{Prop}[Def]{Proposition.}
\newtheorem{Cor}[Def]{Corollary.}
\newtheorem{Lem}[Def]{Lemma.}
\newtheorem{Clm}[Def]{Claim.}
\newtheorem*{Proof}{Proof}
\newtheorem*{Remark}{Remark}
\newtheorem{Ex}[Def]{Example.}

\numberwithin{equation}{section}

\usepackage{secdot}
\sectiondot{section}

\makeatletter
\def\@seccntformat#1{\csname the#1\endcsname.\hs{-.375em}\quad}
\makeatother
\usepackage{tikz}
\usetikzlibrary{matrix}
\begin{document}

\setlength{\abovedisplayskip}{1ex} 
\setlength{\belowdisplayskip}{1ex}

\newcommand{\Uq}{U\hs{-.15em}_q(\mathfrak{g})}
\newcommand{\Ut}{\tilde{U}\hs{-.15em}_q(\mathfrak{g})}
\newcommand{\Up}{U\hs{-.15em}_q^{\hs{.2em}+}
\hs{-.1em}(\mathfrak{g})}
\newcommand{\Um}{U\hs{-.15em}_q^{\hs{.2em}-}
\hs{-.1em}(\mathfrak{g})}
\newcommand{\Uz}{U\hs{-.15em}_q^{{\hs{.2em}0}}
\hs{-.0em}(\mathfrak{g})}

\newcommand{\Utp}{\tilde{U}\hs{-.15em}_q^{\hs{.2em}+}
\hs{-.1em}(\mathfrak{g})}
\newcommand{\Utm}{\tilde{U}\hs{-.15em}_q^{\hs{.2em}-}
\hs{-.1em}(\mathfrak{g})}
\newcommand{\Utz}{\tilde{U}\hs{-.15em}_q^{{\hs{.2em}0}}
\hs{-.0em}(\mathfrak{g})}

\newcommand{\Bb}{\hs{1pt}\overline{\hs{-1pt}B\hs{1pt}}\hs{-1pt}}

\newcommand{\NEWPAGE}{\newpage}
	\renewcommand{\NEWPAGE}{}

\newcommand{\VS}{\vs{4ex}}
	\renewcommand{\VS}{}

\newpage

	\create{center}
		{\Large
			Decomposition Theorem for Product of Fundamental Crystals in Monomial Realization
		}
		
		\vspace{1ex}
                                {\large
                                                Manal Alshuqayr\footnote{ Division of Green Science and Engineering, Sophia University, Kioicho 7-1. Chiyoda-ku, Tokyo 102-8554 Japan, E-mail address: m-alshuqayr@eagle.sophia.ac.jp}} and  {\large Toshiki Nakashima\footnote{ Division of Mathematics, Sophia University, Kioicho 7-1. Chiyoda-ku, Tokyo 102-8554 Japan, E-mail address: toshiki@sophia.ac.jp }

			           }
	\delete{center}

	\thispagestyle{empty}

\begin{abstract}
We consider a product of fundamental crystals in monomial realization. Then we shall show that the product holds crystal structure and describe how it is decomposed into irreducible crystals, which is, in general, different from the decomposition for tensor product of the fundamental crystals. 
\end{abstract}

\section {Introduction }

\TxtBgn
\noindent\indent
        Quantum groups  $\Uq$  are certain families of 
        Hopf algebras that are deformations of universal              
        enveloping algebras of Kac-Moody algebra $\mathfrak{g}$. The theory of crystal 
        bases developed by M. Kashiwara \cite{MK} provides a powerful  
       combinatorial and algebraic tool to study the representations theory of quantum groups.
        Roughly, we can say that crystal base is a basis for certain $\Uq$-module at $q=0$.
        It is shown by Kashiwara that an arbitrary module in the category $\mathcal O_{int}$ \cite{JS} has unique 
       crystal base.
       And a tensor product of such modules also has a crystal base, which means crystal bases hold a natural
       tensor product structure.
\newline 
\noindent\indent
        By Kashiwara operators on  crystal bases, we get the graph structure on a crystal base,
       which is called "crystal graph". By the following simple rule of Kashiwara operators on tensor product (see sect.2):
\begin{align} \label{tenor rule}
\nonumber
       \tilde{e}_{i}\left( b_{1}\otimes \, b_{2} \right)&=\left\{ {\begin{array}{l}
       \tilde{ e}_{i}b_{1}\otimes \, b_{2}\, \, \, \, \, \,\,\mathrm{ if}\,\varphi_{i}\left( 
        b_{1} \right)\ge \varepsilon_{i}\left( b_{2} \right),\, \\ 
        b_{1}\otimes \tilde{\, e}_{i}b_{2}\, \, \, \, \, \, \,\mathrm{ if}\, \varphi 
        _{i}\left( b_{1} \right)<\varepsilon_{i}\left( b_{2} \right), \\ 
\end{array}} \right.\\
\nonumber
       \tilde{f}_{i}\left( b_{1}\otimes \, b_{2} \right)&=\left\{ {\begin{array}{l}
       \tilde{f}_{i}b_{1}\otimes \, b_{2}\, \, \, \, \, \,\,\mathrm{ if}\,\varphi_{i}\left( 
       b_{1} \right)>\varepsilon_{i}\left( b_{2} \right),\, \\ 
       b_{1}\otimes \tilde{\, f}_{i}b_{2}\, \, \, \, \, \, \,\mathrm{ if}\, \varphi 
      _{i}\left( b_{1} \right)\le \varepsilon_{i}\left( b_{2} \right),\\ 
\end{array}} \right.
\end{align}

\noindent
       we know that explicit structure on tensor product of crystal bases.
 
\noindent\indent
       Indeed, a connected component in the crystal graph of a tensor product  corresponds to a simple module
       in the tensor product of modules. 
       That is, let $B(\lambda)$ be a crystal of an irreducible module $V(\lambda)$ $(\lambda $ is dominant),
       then,  for any dominant weight $\lambda,\,\mu$, there exist   dominant weights
       $\lambda_{1},\mathellipsis ,\lambda_{k}$ such that 
\begin{center}     
       $B(\lambda)\otimes B(\mu)\cong B(\lambda_{1}) \oplus \mathellipsis\oplus B(\lambda_{k})$,
\end{center} 
       and
\begin{center}     
       $V(\lambda)\otimes V(\mu)\cong V(\lambda_{1}) \oplus \mathellipsis\oplus V(\lambda_{k})$.
\end{center}

\noindent\indent
      There are several types of realization for crystal bases, $e.g.$, tableaux realization, path realization,
      polyhedral realization, $etc$. Here we will treat the "monomial realization" which is introduced by
      Nakajima \cite{N} and Kashiwara \cite{MK}. Let $\mathcal{M}$ be the set of Laurent monomials in the variables $Y_{i}(n)\, (i \in  I,
              n \in \mathbb{ Z},$ see sect.3): 
\[
       \mathcal{M} \coloneqq\left\{ \prod\limits_{i\in I,n\in \mathbb{ Z}} {Y_{i}\left( n \right)}
      ^{{y}_{i}\left( n 
      \right)} \, ;\, y_{i}(n)\in\mathbb{Z}\, {\hbox{vanish except finitely many }}(i,\, n)\right\} .
\]
      Then we can define the crystal structures on this set of Laurent monomials, and it is shown if a monomial 
      $Y$ satisfies the highest condition, then the connected component in $\mathcal{M}$ including
      $Y$ is isomorphic to the crystal $B\left(wt(Y)\right)$ corresponding to the module $V(wt(Y))$ \cite{MK} as we mentioned above.
      Let us denote this component by $\mathcal{M}(Y)$. Of course, $\mathcal{M}(Y)\otimes \mathcal{M}
      (Y')$ has a crystal structure. But it is not clear whether the set of products
\[
      \mathcal{M}(Y)\cdot\mathcal{M}(Y') \coloneqq\left\{M_{1}\cdot M_{2} \vert\, M_{1}\in \mathcal{M}(Y),
       M_{2}\in \mathcal{M}(Y')\right\},
\]
      holds a crystal structure, where $M_{1}\cdot M_{2}$ means a natural product of Laurent  monomials.

	Let $\mathcal{M}\left(Y_{i}(n) \right) (i \in  I,
              n \in \mathbb{ Z})$ be one of monomial realizations for the fundamental crystal $B(\Lambda_{i})$.

\noindent\indent
      In this article, for type $A_{n}$ we consider the following problems:
\begin{enumerate}
       \item[(1)] \label{crystal structure} dose
       $\mathcal{M}\left(Y_{p}(m) \right)\cdot \mathcal{M}\left(Y_{q}( 1 )\right)$
       $(p,q\in I,\,m\in \mathbb{Z})$ hold a crystal structure?
        \item[(2)] \label{describe the decomposition}
        if the answer for \eqref{crystal structure} is affirmative, describe the decomposition 
       of the crystal 

 $\mathcal{M}\left(Y_{p}(m) \right)\cdot \mathcal{M}\left(Y_{q}( 1 )\right)$.
\end{enumerate} 
       The answer for \eqref{crystal structure} is  affirmative which is shown in 
       Propostion \ref{pos.crystal}, and more general results are obtained in Corollary \ref{general result}; for any $m_1,m_2, \mathellipsis m_l\in  \mathbb{ Z}$,  $p_1,p_2, \mathellipsis p_l
             \in I\, (l>0)$.
                  \begin{center}
             $ \mathcal{M}\left(Y_{p_1}( m_1) \right) \cdot  \mathcal{M}\left(Y_{p_2}(   
              m_2) \right) \cdot  \mathellipsis  \cdot  \mathcal{M}\left(Y_{p_l}( m_l) \right)$ ,
                \end{center}
               holds a crystal structure. Here note that there might be some monomials $M_{1},M_{2}\in \mathcal{M}\left(Y\right)$ and $M'_{1},M'_{2}\in \mathcal{M}
      (Y')$ such that $M_{1}\neq M'_{1}$,  $M_{2}\neq M'_{2}$ and $M_{1}\cdot M_{2}= M'_{1}\cdot M'_{2}$. Therefore, it may happen:
\[
\#\mathcal{M}(Y)\cdot\mathcal{M}(Y')<\#\mathcal{M}(Y)\otimes\mathcal{M}(Y'),
\]
     and then one can deduce that the decomposition for $\mathcal{M}(Y)\cdot\mathcal{M}(Y')$ may differ from the one for $\mathcal{M}(Y)\otimes\mathcal{M}(Y')$. Indeed, the answer for (2) is  given in Theorem \ref{smy.deco.}.
  \begin{enumerate}
                   \item If  $p+q>n$, we have 
                    \[
                   \mathcal{M}\left(Y_{p}(m) \right) \cdot\mathcal{M}\left(Y_{q}( 1 )\right)
                   \cong B(\Lambda_{p}+\Lambda_{q})\,\oplus\bigoplus \limits_{i=\max(p           
                     +q-n,q+1-m)}^{\min(p,q)-1}   
                    B(\Lambda_{p+q-i}+\Lambda_{i})\oplus
                   \mathcal{L}\left(m\ge n-p+2\right)\cdot B(\Lambda_{p+q-n-1}).
                   \]
                
                   \item If $p+q\le n$, we have 
                   \[
                    \mathcal{M}\left(Y_{p}(m) \right)\cdot \mathcal{M}\left(Y_{q}( 1 )\right)
                   \cong B(\Lambda_{p}+\Lambda_{q})\,\oplus\bigoplus \limits_{i=
                     \max(1,q+1-m)}^{\min(p,q)-1}
                   B(\Lambda_{p+q-i}+\Lambda_{i})\oplus
                   \mathcal{L}\left(m\ge q+1\right)\cdot B(\Lambda_{p+q}),
                  \]
                 \end{enumerate}
where $\mathcal{L}(P)=1$ if $P$ is true and $\mathcal{L}(P)=0$ otherwise.
           For example, by this theorem, if $m$ is sufficiently large, we find that 
          $\mathcal{M}(Y_{p}(m))\cdot\mathcal{M}(Y_{q}(1))\cong B(\Lambda_{p})\otimes
          B(\Lambda_{q})$.
           But, if not we find that $\mathcal{M}(Y_{p}(m))\cdot\mathcal{M}(Y_{q}(1))$ is not 
           necessarily isomorphic to $ B(\Lambda_{p})\otimes B(\Lambda_{q})$.
          Indeed, if $m=1$, we know that   $\mathcal{M}(Y_{p}(m))\cdot\mathcal{M}(Y_{q}(1))\cong 
         B(\Lambda_{p}+\Lambda_{q})\varsubsetneq B(\Lambda_{p})\otimes B(\Lambda_{q})$.

       \noindent\indent
      The organization of the article is as follows:
\newline
      In Section 2, we will review on quantum groups, crystal bases, crystals and there
   properties  on tensor product.
      In Section 3,Nakajima's monomial realization of crystals will be introduced , which is our main subject in the article \cite{N}  \cite{MK}.
   In   Section 4 , the explicit form of  monomial realization for the fundamental crystal $B(\Lambda_{i})$ of type $A_{n}$  will be given following \cite{SJD}.
In the last Section, we will state the main results, which are the answers for the questions (1) and (2). In fact, we will show that the product for
      fundamental crystals in monomial realization of type $A_{n}$ holds a crystal structure and describe its decomposition theorem by  classifying the highest weight  
      monomials in the product.

After uploading the article to arXiv, J.Kamnitzer notified us that in their article `` Highest weight for truncated shifted Yangians and product monomial crystals," they treated the same problem for simply-laced cases.
Their method for approaching the problem is some geometric way and then completely different from ours, which is purely combinatorial and direct.
\TxtEnd

{\bf Acknowledgments}: We appreciate M.Nakasuji and Y.Kanakubo for there helpful comments and suggestions.

\section{Preliminaries}
	\noindent\indent
  Let $A=(a_{ij})$ be the Cartan matrix for a simple Lie algebra $\mathfrak{g}$ and let $\mathbb{Q}(q)$ be the rational function field in $q$.

         \DefBgn
        Let $\mathfrak{g}$ be a finite-dimensional simple Lie algebra with a Cartan sub-algebra t, the 
       set of simple roots $\{\alpha_i\in t ^\ast\}_{i\in I}$ and the set of simple coroots 
        $\{h_i\in t\}_{i\in I},$  where $I=\{\,1,\mathellipsis ,n\}\,$. We take an inner product $(\,\,,\,
       \,)$ on $t^\ast$ such that $(\alpha_i,\alpha_i)\in \mathbb{Z}_{>0}$ and 
      $\langle h_i,  \lambda \rangle =2( \alpha_i, \lambda )/(\alpha_i,\alpha_i)$ for 
       $\lambda_i\in t ^\ast $. Let $\{\Lambda_i\}_{i\in I}$ be the dual base of $\{h_i\}$ and set
      $P \coloneqq \sum \mathbb{Z}\Lambda_i$ and  $P^\vee \coloneqq \sum \mathbb{Z}h_i$. 
      Then the $q$-analogue  $U_q(\mathfrak{g})$ is the algebra over  $\mathbb{Q}(q)$ 
      generated by $e_{i},$ $f_{i},$ $k_{i}^{\pm 1}$( $i \in I=\{\,1,\mathellipsis ,n\}\,$) satisfying  
      the relations
\[
		k_i  k_i^{-1} = 1
		= k_i^{-1} k_i, \,\,
		k_i  k_j = k_j  k_i,
	\tag{R1} 
\]
\[
k_i e_j  k_i^{-1} 
		= q_i^{a_{ij}}
		e _j,
								\tag{R2} 
\]
\[
k_i  f_j  k_i^{-1} 
		=  q_i^{a_{ij}}
		f_j,									\tag{R3} 
	\]
\[
		e_i  f_j  - f_j  e_i 
		= \delta_{ij}\,
		\frac{k_i - k_j^{-1}}
		{q_i- q_i^{-1}},		\tag{R4}
	\]

\[
		\begin{array}{l l}
			\displaystyle\sum_{s=0}^{1 - a_{ij}}
			(-1)^s
			\left[\,
			\begin{matrix}
				1 - a_{ij} \\ s
			\end{matrix}
			\,\right]_i
			e_i^{1 - a_{ij} - s}
			e_j
			e_i^{\,s}
			= 0, (i\ne j),
			\tag{R5}
		\end{array}
	\]
	\[
		\begin{array}{l l}
			\displaystyle\sum_{s=0}^{1 - a_{ij}}
			(-1)^s
			\left[\,
			\begin{matrix}
				1 - a_{ij} \\ s
			\end{matrix}
			\,\right]_i
			f_{i}^{1 - a_{ij} - s}
			f_j
			f_{i}^{\,s}
			= 0,  (i\ne j),
			\tag{R6}
		\end{array}
	\]
	where $\delta_{ij}$ is the Kronecker delta.
        \DefEnd

\NtnBgn
	Set now for all $\{\,\alpha_{i}\vert i=1,2,\mathellipsis ,n\}\,$ 
	\create{equation}
		q_i =q^{(\alpha_{i},\alpha_{i})},
	\delete{equation}
	and (for all $a \in \mathbb{Z}$)
	\vs{-2ex}\create{equation}
		[a]_i\coloneqq 
		\frac{q_i^{\,a} - q_i^{-a}}{q_i
		- q_i^{-1}}.
	\delete{equation}
	Define $[n]_{i}!\coloneqq [n]_i[n-1]_i\mathellipsis[2]_i[1]_i$, and
	$
		\left[\,
		\create{matrix}
			a \\ n
		\delete{matrix}
		\,\right]_i \coloneqq \frac{[a]_{i}!}{[n]_{i}![a-n]_{i}!}\,(a,n\in \mathbb{Z}_{\ge 0},a\ge n) .
	$
\NtnEnd

	\iftrue
 \TxtBgn
      Let $M$ be a  
      finite dimensional $U_q(\mathfrak{g})$-module and for $\lambda\in P$ set $M_\lambda=\{v\in M  \vert q^{h}
      v=q^{\langle h, \lambda\rangle }(\forall  h\in P^{\vee})\}$, which is called a weight space 
      of weight $\lambda$. Then we know that  $M=\bigoplus_{\lambda \in P}M_{\lambda } $.
\TxtEnd

\DefBgn 
                  Let $M=\bigoplus_{\lambda \in P}M_{\lambda } $ be a 
                 finite dimensional  $\Uq$-module.
                  For each $i\, \in \, I$, every weight vector $u\, \in \, M_{\lambda }\, 
                  (\lambda \in \, \, wt(M))$ is written in the form 
                        \begin{align} \label{form of u}
   u=u_{0}\, +\, f_{i}u_{1}\, + \mathellipsis+\, f_{i}^{(N)}u_{N}, 
                       \end{align}
               where $N \in \mathbb{Z} _{\ge 0 }$ and $ u_{k}\in M_{\lambda +k\alpha 
                  _{i}}\cap\Ker e_{i}$ for any $k=0,1,\mathellipsis ,N$.
                  Here, each $u_{k}$ in the expression is uniquely determined by $u$, and 
                  $u_{k}\ne 0$ only if $\lambda (h_{i})+ k>0$. 
                The\itshape{ Kashiwara operators}  $\tilde{e}_{i},\tilde{f}_{i}\in End_{\mathbb{Q}(q)}(M)\left( i\in I \right)$ are defined  for $u=\sum\limits_{k=0}^{N}f_{i}^{(k)}u_{k}\in M_{\lambda}$ by                       
                          \begin{align}
                  \tilde{e}_{i}u=\sum\limits_{k=1}^N f_{i}^{(k-1)} \, u_{k},\, \, \, \, 
                  \, \, \tilde{f}_{i}u=\sum\limits_{k=0}^N f_{i}^{(k+1)} \, u_{k}\, .
                         \end{align}
\DefEnd

 \TxtBgn
             Let $A$ be the subring of $\mathbb{Q}(q)$ defined by $A=\{f(q)/g(q)\,\vert\, f(q),g(q)
             \in \mathbb{Q}[q], g(q)\neq 0\}$.
\TxtEnd

\DefBgn 
                 Let $M$ be a finite dimensional  $\Uq$-module and $L$ be a 
                 free $A$-submodule. $ L $ is called a \itshape{crystal lattice} if 
           \vs{-1.5ex}     \begin{enumerate}
                \item $L$ generates $M$ as a vector space over $\mathbb{Q}(q)$,
                 \item $L=\bigoplus_{\lambda \in P}L_{\lambda },$ where 
                 $L_{\lambda }=L \cap M_{\lambda }$ for all $\lambda \in P$,
                  \item $\tilde{e}_{i}L\subset L,\, \tilde{f}_{i}L\subset L$
                   for all $ i\in I .$
                 \end{enumerate}
\DefEnd

\vs{-1ex}
    \DefBgn \label{crystal basis}
               A \itshape{crystal base} of a  finite dimensional  $\Uq$-module $M$  is a    
                pair $(L, B)$ such that 
           \vs{-1.5ex}  \begin{enumerate}
              \item $L$ is a crystal lattice of $M$,
              \item $B$ is an $\mathbb{Q}$-basis of $L/qL$,
               \item $B=\bigsqcup_{\lambda \in P }B_{\lambda},$
                 where $ B_{\lambda }= B\cap ( L_{\lambda }/q L_{\lambda }),$ 
                 \item $\tilde{e}_{i}\,B\subset B\cup \left\{ 0 \right\},\,
                 \tilde{f}_{i}B\subset B\cup \left\{ 0 \right\}$ for all $ i\in I,$
                 \item  for any $ b,\, b' \in B$ and $ i\in I,$ we have $\tilde{f}_{i}b=b'$
                 if and only if $ b=\tilde{e}_{i}b'$.
                  \end{enumerate}
\DefEnd

\vs{-.5ex}
\TxtBgn 
          Let $P_{+}\coloneqq \{\sum\limits_{i=1}^{n}m_{i}\Lambda_{i}\, \vert\, m_{i}\in \mathbb{Z} 
           _{\ge 0 }
          (\forall i\in I)\}$ be the set of dominant weights. For any finite dimensional irreducible 
          $\Uq$-module $L$, there exists a unique dominant weight $\lambda\in P_{+}$ such 
         that $L \cong V(\lambda)$, where $V(\lambda)\coloneqq \Uq / \sum\limits_{i\in I} \Uq e_{i}
         +\sum\limits_{i\in I} \Uq(t_{i}-q_{i}^{\langle h_i,  \lambda \rangle})+\sum\limits_{i\in I} \Uq 
          f_{i}^{1+\langle h_i,  \lambda \rangle}$.

        Here let $\pi _{\lambda}\coloneqq \Uq \rightarrow V(\lambda)$ be the natural projection and 
         set $u_{\lambda}\coloneqq \pi _{\lambda}(1)$, which is the highest weight vector of 
          $ V(\lambda)$.

         Let us define for $\lambda\in P_{+}$, 
            \begin{align}
           \nonumber
         L(\lambda)\coloneqq &\sum\limits_{\substack{ i_{1},
       \mathellipsis ,i_{l}\in I\\ l\ge 0}}A \tilde{f}_{i_1}\mathellipsis\tilde{f}_{i_l}u_{\lambda},\\
           \nonumber
        B(\lambda)\coloneqq &\left\{\tilde{f}_{i_1}\mathellipsis\tilde{f}_{i_l}u_{\lambda}{\rm mod}\, qL
        (\lambda)\,\vert\, i_{1},
       \mathellipsis ,i_{l}\in I, l\ge 0\right\}\setminus\{0\}.
  \end{align}
         \TxtEnd

\vs{-1.5ex}
\ThmBgn  \cite{M1,M2}
       The pair $(L(\lambda),B(\lambda))$ is a crystal base of $V(\lambda)$.
\ThmEnd

\ThmBgn\cite{M1,M2}
          Let $M_{j}$ be a  finite dimensional $\Uq$-module and let $\left( L_{j},
          B_{j} \right)\, \, $be a crystal basis of $M_{j} \, (j\, =\, 1,\, 2)$. Set 
          $L={L_{1}\otimes }_{A_{0}}L_{2}$ and $B=\, B_{1}\times B_{2}$ .

        Then $\left( L, B \right)$ is a crystal basis of ${M_{1 }\otimes 
        }_{A_{0}}M_{2 }$, where the action of Kashiwara operators$\, 
        \tilde{e}_{i}$ and $\tilde{ f}_{i}$ on $B \, (i \in I)$ is given 
        by 
         \begin{align} \label{tenor rule}
         \nonumber
        \tilde{e}_{i}\left( b_{1}\otimes \, b_{2} \right)&=\left\{ {\begin{array}{rcl}
         \tilde{ e}_{i}b_{1}\otimes \, b_{2}&\mbox{ if}&\varphi_{i}\left( 
         b_{1} \right)\ge \varepsilon_{i}\left( b_{2} \right),\\ 
         b_{1}\otimes \tilde{e}_{i}b_{2}&\mbox{ if}& \varphi 
        _{i}\left( b_{1} \right)<\varepsilon_{i}\left( b_{2} \right), 
            \end{array}} \right.\\
         \nonumber
        \tilde{f}_{i}\left( b_{1}\otimes \, b_{2} \right)&=\left\{ {\begin{array}{rcl}
         \tilde{f}_{i}b_{1}\otimes \, b_{2}&\mbox{ if}&\varphi_{i}\left( 
         b_{1} \right)>\varepsilon_{i}\left( b_{2} \right),\\ 
         b_{1}\otimes \tilde{f}_{i}b_{2}&\mbox{ if}&\varphi 
       _{i}\left( b_{1} \right)\le \varepsilon_{i}\left( b_{2} \right), 
         \end{array}} \right.
           \end{align}
      and we have
         \begin{align}
           \nonumber
         wt\left( b_{1}\otimes \, b_{2} \right)&=\, wt\left( b_{1} \right)+\, wt(\, 
         b_{2}),\\
         \nonumber
        \varepsilon_{i}\left( b_{1}\otimes \, b_{2} \right)&=\max (\varepsilon 
       _{i}\left( b_{1} \right){,\varepsilon }_{i}\left( b_{2} 
        \right)-\langle\hs{1pt} \ h_{i}, wt\left( b_{1} \right)\rangle),\\
           \nonumber
        \varphi_{i}\left( b_{1}\otimes \, b_{2} \right)&=\max (\varphi_{i}\left( 
        b_{1} \right),\varphi_{i}\left( b_{2} 
        \right)-\langle\hs{1pt} \ h_{i}, wt\left( b_{2} \right)\rangle).
              \end{align}
            \begin{center}
        Here, we write $b_{1}\otimes \, b_{2}$ for $(b_{1},\, b_{2})\in B_{1}\times 
        B_{2}$ and we understand $b_{1}\otimes 0 =0\otimes  b_{2} =0$.
            \end{center}
\ThmEnd

\DefBgn \cite{M}
Let $I$ be a finite index set and let $A=(a_{ij})_{i,j\in I}$ be a Cartan matrix and $P$ a corresponding weight lattice.   \itshape{Crystal} associated with the Cartan matrix $A$  is a set $B$ together with the maps $wt:B\to P$, $  \tilde{e}_{i},\,  \tilde{f}_{i}:B\to B\cup \{0\}$, and $\varepsilon _{i},\,\varphi_{i}:B\to\mathbb{Z}\cup\{-\infty\}(i\in I)$ satisfying the following properties:
\begin{enumerate}
\item $\varphi_{i}(b)=\varepsilon_{i}(b)+\langle h_{i},\,wt(b)\rangle$ for all $i\in I$,
\item $wt(  \tilde{e}_{i}b)=wt\, b+\alpha_{i}$ if $ \tilde{e}_{i}b\in B$,
\item $wt(  \tilde{f}_{i}b)=wt\, b-\alpha_{i}$ if $ \tilde{f}_{i}b\in B$,
\item $\varepsilon_{i}(  \tilde{e}_{i}b)=\varepsilon_{i} ( b)-1,\,\varphi_{i}(  \tilde{e}_{i}b)=\varphi_{i}( b)+1$ if $\tilde{e}_{i}b\in B$,
\item  $\varepsilon_{i}(  \tilde{f}_{i}b)=\varepsilon_{i} ( b)+1,\,\varphi_{i}(  \tilde{f}_{i}b)=\varphi_{i}( b)-1$ if $\tilde{f}_{i}b\in B$,
\item $ \tilde{f}_{i}b=b'$ if and only if $b=\tilde{e}_{i}b'$ for $b,\,b'\in B$, $i\in I$,
\item if $\varphi_{i}(b)=-\infty$ for $b\in B$, then $  \tilde{e}_{i}b=  \tilde{f}_{i}b=0$.
\end{enumerate}
\DefEnd

\DefBgn  
             Take $B$ as the set of vertices and the $I$-colored 
             arrows on $B$ by
           \begin{center}
              $
             b\buildrel i \over \longrightarrow b'\,$ if and only if
              $\tilde{f}_{i}b =b^{'}(i\in I).$
                \end{center}
               Then $B$ is given an $I$-colored oriented graph structure called the 
                 \itshape{crystal graphs} of $M$.
\DefEnd

                  \section{Nakajima's monomials }
               \TxtBgn
\noindent\indent
In this section, we recall the crystal structure on the set of monomials discovered by H. Nakajima \cite{N}. Our exposition follows that of M.Kashiwara \cite{MK}.

                     \TxtBgn
              Let $\mathcal{M}$ be the set of Laurent monomials in the variables $Y_{i}(n)\, (i \in  I,
              n \in \mathbb{ Z}):$ 
                       \[
       \mathcal{M} \coloneqq\left\{ \prod\limits_{i\in I,n\in \mathbb{ Z}} {Y_{i}\left( n \right)}
      ^{{y}_{i}\left( n 
      \right)} \, ;\, y_{i}(n)\in\mathbb{Z}\, {\hbox{vanish except finitely many }}(i,\, n)\right\} .
\]
           
\noindent\indent
            We shall define a crystal structure of crystal on  $\mathcal{M}$. Let $c=
                      {(c_{ij})}_{i\ne j\in I}$, be a set of integers such that 
           \[
            c_{ij}+\, c_{ji} = 1.
             \]
          We set 
       \[
            A_{i}(n)\, =\, Y_{i}(n)Y_{i}(n\, +\, 1)\, \prod\limits_{j\ne i} {Y_{j}(n\, 
            +c_{ji})}^{\langle h_{j},\alpha_{i}\rangle\,}.
      \]
            For a monomial  $M\, =\, \prod\limits_{i\,\in I,n\,\in\,\mathbb{ Z}} {Y_{i}\left( n 
               \right)}^{y_{i}\left( n \right)}\in\mathcal{M} $, we set
            \begin{align}
              \nonumber
            wt\left( M \right)&=\, \sum\nolimits_i ( \sum\nolimits_n {Y_{i}(n))} \Lambda 
             _{i}\, \, ,\\
               \nonumber
             \varphi_{i}\left( M \right)&=\, max\,\{\, \sum\nolimits_{k\le n} {y_{i}(k)} ;\, 
              n\in \,\mathbb{ Z}\},\\
              \nonumber
             \varepsilon_{i}\left( M \right)&=\, max\{\,-\sum\nolimits_{k> n} {y_{i}(k)} ;\, 
              n\in \, \mathbb{ Z}\}.
                  \end{align} 
                       We define 
                           \begin{align}
                        \nonumber
                   \tilde{ f}_{i}\left( M \right)&=\left\{ {\begin{array}{rcl}
                  0\quad\quad\quad&\mbox{ If}& \varphi_{i}\left( M \right)=0,\\ 
                   A_{i}{(n_{f})}^{-1}M&\mbox{  If}&\varphi 
                  _{i}\left( M \right)>0,
                          \end{array}} \right.\\
                         \nonumber
                         \tilde{ e}_{i}\left( M \right)&=\left\{ {\begin{array}{rcl}
                      0\quad\quad&\mbox{ If}& \varepsilon_{i}\left( M \right)=0, \\ 
                 A_{i}(n_{e})M&\mbox{  If}&\varepsilon 
                 _{i}\left( M \right)>0, 
                               \end{array}} \right.
                           \end{align}
                            where 
                             \begin{align}
                           \nonumber
                        n_{f}&=min\left\{n;\, \varphi_{i}\left( M 
                \right)\mathrm{=}\sum\nolimits_{k\le n} {y_{i}(k)}\right\} \, ,\\
                          \nonumber
                       &=min\left\{n;\, \varepsilon_{i}\left( M 
                         \right)\mathrm{=-}\sum\nolimits_{k>n} {y_{i}(k)}\right\} \, ,\\
                \nonumber
                        n_{e}&=max\left\{ \mathrm{n;\, }\varphi_{i}\left( M 
                     \right)\mathrm{=}\sum\nolimits_{k\le n} {y_{i}\left( k \right)} \right\},\\     
                        \nonumber
                  &=max\left\{n;\, \varepsilon_{i}\left( M 
                         \right)\mathrm{=-}\sum\nolimits_{k>n} {y_{i}(k)} \right\}\,.
                 \end{align}

                           Note that $y_{i}(n_{f}) > 0,\,y_{i}(n_{f}+1)\le 0$ and
                      $  y_{i}(n_{e}+1)< 0,\,y_{i}(n_{e})\ge 0$. 
 
                    Let us denote by  $\mathcal{M}_{c}$ the crystal $\mathcal{M}$ associated with c. 
                           \TxtEnd

\TxtBgn
   We shall 
             call $I$-crystal when we want to emphasize the index set $I$ of simple 
                  roots.

\noindent\indent
           For an $I$-crystal $B\, $and a subset $J$ of $ I$,  let us denote by $\Psi_{J} (B)$ the
            $J$-crystal $B$ with the induced crystal structure from $B$.

\noindent\indent
            A crystal $B$ is called \itshape{ semi-normal} if for each $i\in I$  the 
           $\{\,i\}\,$-crystal $\Psi 
           _{\{\,i\}\,}(B)$  is a crystal associated with a finite dimensional  module. This is 
           equivalent to saying that $\varepsilon_{i}(b)=max \{\,n \in \, N\, ;\, 
           \tilde{e}_{i}^{n}b \ne 0\}$ and $\varphi_{i}(b)=max\{\,n \in N\, 
           ;\, \tilde{f}_{i}^{n}b\ne 0\}$ for any $b\, \in \, B$ and $i\, \in I$.
         \TxtEnd

\noindent\indent
              A crystal $ B$ is called \itshape{ normal} if for any subset $J$  of  $I$ such that 
              $\{\,\alpha_{i}\, ;\, i\in J\}\,$ is a set of simple roots for finite-dimensional 
               semisimple Lie algebra $\mathfrak{g_{\mathrm{J}}}$, $\Psi_{J}(B)$ is the crystal associated with 
             a finite dimensional  $U_q(\mathfrak{g_{\mathrm{J}}})$-module, where $U_q(\mathfrak
                {g_{\mathrm{J}}})$ is 
               the associated quantum group with $\{\,\alpha_{i}\, ;\, i\, \in  \, J\}\,$.

                  \PrpBgn
                      $\mathcal{M}_{c}$ is a semi-normal crystal. 
               \PrpEnd

            \TxtBgn
                     For a family of integers ${(m_{i})}_{i\in I}$, let us set $c'\, =\, 
                   {({c'}_{ij})}_{i\ne j\in I}$, by ${c'}_{ij}=\, c_{ij}+m_{i}\, -m_{j}$. Then 
                   the map $Y_{i}\left( n \right)\longmapsto \, Y_{i}\left( n+m_{i} \right)$ 
                   gives a crystal isomorphism  $\mathcal{M}_{c}\buildrel \sim\over\longrightarrow 
                  \mathcal{M}_{c'}$.
                  Hence if the Dynkin diagram 
                  has no loop, the isomorphism class of the crystal $\mathcal{M}_{c}$ does not 
                 depend on the choice of $c$. 
             \TxtEnd

             \ThmBgn \cite{MK}
                 For a highest weight vector $M \in \mathcal{M}_{c}$, the 
                  connected component of  $\mathcal{M}_{c}$  containing $M$ is isomorphic to
                 $B(wt(M))$. 
               \ThmEnd

                     \CorBgn  \cite{MK}  \label{normality} 
$\mathcal{M}_{c}$  is a normal crystal.
                    \CorEnd

                     \PrpBgn \label{M crystal }\begin{enumerate}
                    \item For each $i \in  I$, $\mathcal{M}$  is isomorphic to a crystal graph of a
                      $U_{q}(\mathfrak{sl}_2) $-module.
                    \item Let $M$ be a monomial with weight $\lambda $, such that $\tilde{ e}_{i}M= 0$
                    for all $ i \in  I$, and let $\mathcal{M}(\lambda )\, $be the connected component of
                     $\mathcal{M}$ containing $M$. Then there exists a crystal isomorphism
                       \end{enumerate}
                       \begin{center}
                 $\mathcal{M}\left( \lambda \right) \buildrel\sim\over\longrightarrow\,
                   B\left( \lambda \right)$\, \, \, \, given by
                  \, \,$ M\longmapsto v_{\lambda }$.
                    \end{center}
                  \PrpEnd

\RmkBgn
Not that, the condition $\tilde{ e}_{i}M= 0$
                    for all $ i \in  I$, is called {\it the highest weight} condition, which is, equivalent to the condition $\varphi_{i}\left( M \right)=0$   for all $ i \in  I$.
\RmkEnd

               \section{Monomial Realization of the Fundamental Crystals of Type $A_{n}$}
                 \TxtBgn
\noindent\indent
                       We refer to \cite {SJD}.
                        Let $I\, =\, 
                       \{1,\,\mathellipsis ,\, n\}\,$ and let $A\, =\, {(a_{ij})}_{i,j\in I}$ be the 
                        Cartan matrix of type $A_{n}$. Here, the entries of $A$ are given by 
            \[
                         a_{ij}=\left\{ {\begin{array}{rcl}
                         2&\mbox{if}& i=j, \\ 
                        -1&\mbox{if}&\left| i-j \right|=1, \\ 
                        0&  &\mbox{otherwise}.
                                                  \end{array}} \right.
               \]

                   Let $U_{q}(\mathfrak g)\, =\, U_{q}({sl}_{n+1})$ the
                      corresponding quantum 
                      group. For simplicity, we take the set $C\, ={(c_{ij})}_{i\ne j}$ to be 
                       \begin{center}
                    $c_{ij}=\left\{ {\begin{array}{rcl}
                    0& \mbox{if}& i>j , \\ 
                      1 &\mbox{if}& i<j,
                       \end{array}} \right.$
                        \end{center}
                     and set ${Y_{0}(m)}^{\pm 1}={Y_{n+1}(m)}^{\pm 1}
                       =1$ for all $ m\in  \mathbb{ Z}$. Then for $i\ \in  I$
                     and $m\in  \mathbb{ Z}$, we have
                    \begin{center}
                    $A_{i}\, (m)\, =Y_{i}(m)Y_{i}(m\, +\, 1){Y_{i-1}(m\, +\, 
                   1)}^{-1}{Y_{i+1}(m)}^{-1}$.
                      \end{center}

                   To characterize $\mathcal{M} \left( \lambda \right)$, we first 
                   focus on the case when $\lambda =\, \mathrm{\, }\Lambda_{k}$. Let $M_{0} 
                  =Y_{k}\left( m \right)$  for $ {m\in  \mathbb{ Z}}$.
                Hence $\tilde{e}_{i}M_{0}=0$ for all $i \in I $ 
                and the connected component containing $M_{0}$ is isomorphic to $B(\Lambda 
              _{k})$ over $U_{q}(\mathfrak g).$ For simplicity, we will take $M_{0}=
                  Y_{k}(1)$, even if that does not make much difference.

                 \PrpBgn\label{factor}
                 For $k = 1,\mathellipsis ,n$, let $M_{0}=Y_{k}(1)$ be the monomial of weight 
                 $\Lambda_{k}$ such that $\tilde{ e}_{i}M_{0} =0$  for all
                   $ i\in I$. Then the connected component $\mathcal{ M}(M_{0})=\mathcal{ M}(Y_{k}(1))$ of $M$ 
                  containing $M_{0}$ is characterized as 
           \[
                 \mathcal{M} \left( (Y_{k}(1) \right)=
                 \left\{ \prod\limits_{j=1}^r {Y_{a_{j}}\left( m_{j-1}
                 \right)^{-1}Y_{b_{j}}\left( m_{j} \right)}\,\,\vert \left. {\begin{array}{l}
                 0\le a_{1}<b_{1}<a_{2}<\mathellipsis <a_{r}<b_{r}\le n+1,\\ 
                 k+1=m_{0}>m_{1}>\mathellipsis >m_{r-1}>m_{r}=1 ,\\ 
                a_{j}+m_{j-1}=b_{j}+m_{j}\, \, \, \forall j=1,\mathellipsis ,r\le k. \\ 
                                 \end{array}} \right\} \right.
              \]
            \PrpEnd

                  \RmkBgn  
                       If we take $M_{0}=Y_{k}(N)$, then we have only to 
                       modify the condition for $m_{j}$'s as follows: 
                             \begin{center} 
                       $k +N=\, m_{0}>m_{1}>\mathellipsis \, >m_{r-1}\, >m_{r}=\, N$.
                              \end{center}
                         $For \,\,i\, \in \, I$  and  ${m\in \mathbb{ Z}}$, we introduce
                                                              new variables
                              \begin{center}
                             $ X_{i}\, (m)\, =\, {Y_{i-1}(m\, +\, 1)}^{-1}Y_{i}(m).\, $
                                   \end{center}
                               Note that $X_{n+1}(m)=\frac {1}{Y_{n}(m+1)}$.
                          \vs{.5ex}
                                  Using this 
                                   notation, every monomial
                               \begin{center}
                            $M\, =\prod\nolimits_{j=1}^r {\, {Y_{a_{j}}\left( 
                           m_{j-1} \right)}^{-1}Y_{b_{j}}\left( m_{j} \right)\in M\left( Y_{k}(1) 
                              \right)} $,
                              \end{center}
                        may be written as
                              \begin{center}
                       $\, M=\prod\limits_{j=1}^r \, X_{a_{j}+1}\left( m_{j-1}-1 
                       \right)X_{a_{j}+2}\left( m_{j-1}-2 \right)\mathellipsis X_{b_{j}}\left( 
                        m_{j} \right)$.
                                 \end{center}
                      For example, we have $Y_{k}\left( N
                       \right)=X_{1}\left(  k+N-1  \right)X_{2}\left(k+N-2  \right)\mathellipsis 
                    X_{k}\left( N\right)$.
                     Now, it is straightforward to verify that we have another characterization 
                   of the crystal$\,\mathcal{ M}\left(Y_{k}(1)\right)$ .
                      \RmkEnd

               \CorBgn \label{form of M}  
                  For $k\, =\, 1,\mathellipsis n$ we have 
         \[
                 \mathcal{M}(Y_{k}(N)) =\left\{\, X_{i_{1}}( k+N-1 
                  )X_{i_{2}}( k+N-2 )\mathellipsis \, X_{i_{k}}( N
                    )\vert  \,\, 1\le i_{1}\, <i_{2}\, <\mathellipsis \, <\, i_{k}\le \, 
                       n+1\right\}.
           \]
                     \CorEnd

\ExpBgn \label{crystal graph}
              Let $\mathfrak{g}=A_{5}$, and          $c_{ij}=\left\{ {\begin{array}{rcl}
                    0& \mbox{if}& i>j , \\ 
                      1 &\mbox{if}& i<j.
                       \end{array}} \right.$

 \begin{enumerate}
                     \item  The crystal    $B(\Lambda_{5})\cong\mathcal{M}\left( Y_{5}(m) \right)$ is given 
                  as follows.
\[
                \mathcal{M}\left( Y_{5}(m) \right)=\left\{\, X_{i_{1}}\left( 4+m
                  \right)X_{i_{2}}\left( 3+m\right)X_{i_{3}}\left( 2+m\right)
                  X_{i_{4}}\left( 1+m\right)X_{i_{5}}\left( m\right)
                 \vert \,1\le i_{1}\, <i_{2} <i_{3} <i_{4} <i_{5}\le 6 \right\},
\]
\begin{center}
                \begin{tikzpicture}
        \matrix (m) [matrix of math nodes,row sep=3.5em,column sep=1.5em,minimum width=2em]
  {
      Y_{5}(m)&Y_{5}(m+1)^{-1}Y_{4}(m+1)&Y_{4}(m+2)^{-1}Y_{3}(m+2)&Y_{3}(m+3)^{-1}Y_{2}   (m+3)\\ &  & Y_{1}(m+5)^{-1}&Y_{2}(m+4)^{-1}Y_{1}(m+4)
       \\};
         \path[-stealth]
          (m-1-1)   edge  node  [above]{$\tilde{f}_{5}$}(m-1-2)
          (m-1-2)  edge  node [above] {$\tilde{f}_{4}$} (m-1-3)
          (m-1-3) edge node [above]{$\tilde{f}_{3}$} (m-1-4)
          (m-1-4) edge node [right] {$\tilde{f}_{2}$} (m-2-4)
          (m-2-4) edge  node  [above] {$\tilde{f}_{1}$} (m-2-3)   ;
\end{tikzpicture}
\end{center}
 \[\text{
			{\bf Figure 3}
		}\]
\newpage
    \item      The crystal    $B(\Lambda_{2})\cong\mathcal{M}\left( Y_{2}(1) \right)$ is given 
                  as follows.

               $\mathcal{M}\left( Y_{2}(1) \right)=\left\{\, X_{i_{1}}\left( 2 
                  \right)X_{i_{2}}\left( 1 \right)\vert \,1\le i_{1}\, <i_{2}\le 6 \right\}$,
\begin{center}
 \begin{tikzpicture}
        \matrix (m) [matrix of math nodes,row sep=1.5em,column sep=1.5em,minimum width=1em]
  {
     & Y_{2}(1)\\ &Y_{3}(1)^{-1}Y_{2}(2)^{-1}Y_{1}(2)\\Y_{4}(1)Y_{3}(2)^{-1}Y_{1}(2) & &Y_{3}(1)Y_{1}(3)^
{-1}\\ Y_{5}(1)Y_{4}(2)^{-1}Y_{1}(2)& &Y_{4}(1)Y_{3}(2)^
{-1}Y_{2}(2)Y_{1}(3)^{-1}\\Y_{5}(2)^{-1}Y_{1}(2)
&Y_{5}(1)Y_{4}(2)^{-1}Y_{2}(2)Y_{1}(3)^{-1}&Y_{4}(1)Y_{2}(3)^{-1}
\\Y_{5}(2)^{-1}Y_{2}(2)Y_{1}(3)^{-1}& &Y_{5}(1)Y_{4}(2)^{-1}Y_{3}(2)Y_{2}(3)^{-1}
\\ Y_{5}(2)^{-1}Y_{3}(2)Y_{2}(3)^{-1}&  &Y_{5}(1)Y_{3}(3)^{-1}
\\ &Y_{5}(2)^{-1}Y_{4}(2)Y_{3}(3)^{-1}& \\  &Y_{4}(3)^{-1}&  
       \\};
         \path[-stealth]
          (m-1-2)   edge  node [right] {$\tilde{f}_{2}$}(m-2-2)
          (m-2-2)  edge  node[above] {$\tilde{f}_{3}$} (m-3-1)
                        edge node[above] {$\tilde{f}_{1}$}(m-3-3)
          (m-3-1) edge node[left]{$\tilde{f}_{4}$} (m-4-1)
                     edge node[above] {$\tilde{f}_{1}$}(m-4-3)
            (m-3-3) edge node[right] {$\tilde{f}_{3}$} (m-4-3)
           (m-4-1)edge node[left]{$\tilde{f}_{5}$} (m-5-1)
                       edge node[above]{$\tilde{f}_{1}$}(m-5-2)
            (m-4-3) edge node[above]{$\tilde{f}_{4}$}(m-5-2)
                          edge node[right]{$\tilde{f}_{2}$}(m-5-3)
           (m-5-1) edge node [left]{$\tilde{f}_{1}$}(m-6-1)
           (m-5-2) edge node [above]{$\tilde{f}_{2}$}(m-6-3)
                        edge node[above]{$\tilde{f}_{5}$}(m-6-1)
           (m-5-3) edge node[right]{$\tilde{f}_{4}$}(m-6-3)
           (m-6-3) edge node [above]{$\tilde{f}_{5}$}(m-7-1)
           (m-6-1) edge node [left]{$\tilde{f}_{2}$}(m-7-1)
           (m-6-3) edge node [right]{$\tilde{f}_{3}$}(m-7-3)
           (m-7-1) edge node[below]{$\tilde{f}_{3}$}(m-8-2)
           (m-7-3) edge node[below]{$\tilde{f}_{5}$}(m-8-2)
           (m-8-2) edge node[left]{$\tilde{f}_{4}$}(m-9-2)
 ;
\end{tikzpicture}
\end{center}
 \[\text{
			{\bf Figure 4}
		}\]
\end{enumerate}
\ExpEnd

                  \section{Product of Fundamental Crystals in Monomial Realization}
                   \TxtBgn
\noindent\indent
                     In this section we study the crystal structures of product for fundamental 
                     crystals in monomial realization of type $A_{n}$. We will show that such a
                      product has a crystal structure and  we introduce the
                      decomposition theorem for product of fundamental crystals in monomial realization. 
                    \TxtEnd

              \subsection{The Product  $ \mathcal{M}\left(Y_{p}( 
                    m) \right)\cdot \mathcal{M}\left(Y_{q}( 1 )\right)$ }

                   \TxtBgn
\noindent\indent
                       We use the same notation as the previous sections.
                     \TxtEnd

             \LemBgn \label{crystal of product}
                    Let $M_{1}\in\mathcal{M}\left(Y_{p}( m )\right)$, 
                    $M_{2}\in\mathcal{M}\left(Y_{q}( 1) \right)$ with $m\ge 1$, and  $p,q \in\{1,\mathellipsis ,n\}\,$.
                     Then for any $i\in I$, we have
           \[
                 \tilde{e}_{i}\, \left( M_{1}\cdot{M}_{2} \right)\in  \mathcal{M}(Y_{p}\left( m \right))
                \cdot \mathcal{M}(Y_{q}\left( 1 \right))\cup \, \{0\},
             \]
          \[
               \tilde{f}_{i}\, \left( M_{1}\cdot{M}_{2} \right)\in \mathcal{M}(Y_{p}\left( m \right)) 
              \cdot \mathcal{M}(Y_{q}\left( 1 \right))\cup \, \left\{ 0 \right\}.
             \]
          \LemEnd

         \PrfBgn
                    By the explicit description in Corollary  \ref{form of M}, we can write 
       \begin{align}
                \nonumber
                     M_{1}&=X_{i_{1}}\left( p+m-1 \right)X_{i_{2}}\left( p+m-2 
                     \right)\mathellipsis X_{i_{p}}(m),\\
              \nonumber
                    M_{2}&=X_{j_{1}}\left( q \right)X_{j_{2}}\left( q-1 \right)\mathellipsis 
                    X_{j_{q}}\left( 1 \right).
            \end{align}
                  The factor $Y_{i}\left(l\right)^{-1}(l\in\mathbb{Z})$ can appear at most once in
                  $M_{1}$ and $M_{2}$ by Proposition \ref{factor} respectively
                             then we obtain 

    \[
                    \tilde{e}_{i}\, \left( M_{1}\cdot{M}_{2} \right)=\left\{ {\begin{array}{l}
                    \left( \tilde{e}_{i}M_{1} \right)\cdot M_{2}, \\     
                    M_{1}\cdot \left( \tilde{e}_{i}M_{2} \right), \\ 
                     0. \\ 
                 \end{array}} \right.
           \]
                  Indeed, since if both  $Y_{i}\left( l_{1} \right)^{-1}$ appears in $M_{1}$
                     and $Y_{i}\left( l_{2} \right)^{-1}$ appears in  $M_{2}$, by the definition 
                      of $ \tilde{e}_{i}$ we have
             \begin{align}
                  \nonumber
                    \tilde{e}_{i}\, \left( M_{1}\cdot{M}_{2} \right)& =A_{i}\left( l_{1} \right) 
                     \left( M_{1}\cdot{M}_{2} \right)\,\,\,\mathrm{or}\,\,\,A_{i}\left( l_{2} \right)
                     \left( M_{1}\cdot{M}_{2} \right)\\
                  \nonumber
                    &= \left( \tilde{e}_{i}M_{1} \right)\cdot M_{2}
                     \,\,\,\mathrm{or}\,\,\,  M_{1}\cdot \left( \tilde{e}_{i}M_{2} \right).
                 \end{align}
                       If one of them appears in $M_{1}$ or $M_{2}$,
                                   we also have 
               \begin{center}
                  $\tilde{e}_{i}\, \left( M_{1}\cdot{M}_{2} \right)=
                       \left( \tilde{e}_{i}M_{1} \right)\cdot M_{2}$
                        or $  M_{1}\cdot \left( \tilde{e}_{i}M_{2} \right)$ .
                \end{center}

                    Similarly, by Proposition \ref{factor} we have 
                \[
                        \tilde{f}_{i}\, \left( M_{1}\cdot{M}_{2} \right)=\left\{ {\begin{array}{l}
                        \left( \tilde{f}_{i}M_{1} \right)\cdot M_{2}, \\ 
                         M_{1}\cdot \left( \tilde{f}_{i}M_{2} \right), \\ 
                                     0. \\ 
                           \end{array}} \right.
                 \]
                \PrfEnd

                  \PrpBgn \label{pos.crystal}
             The product $ \mathcal{M}\left(Y_{p}(  m) \right)\cdot \mathcal{M}\left(Y_{q}( 1 )\right)$
             possesses the crystal structure and it is decomposed into a direct sum of crystals, that is, 
            there exist dominant integral weights $\lambda_{1},\lambda_{2},\mathellipsis,\lambda_{k}\in
            P_{+}$ \,such that 
          \begin{center}
           $ \mathcal{M}\left(Y_{p}( 
                    m) \right) \cdot\mathcal{M}\left(Y_{q}( 1 )\right)
           \cong B\left( \lambda_{1} \right)\oplus B\left( \lambda_{2} \right) \oplus\mathellipsis
           \oplus  B \left( \lambda_{k} \right)$.
           \end{center}
              \PrpEnd

           \PrfBgn
           The former half of the statement is clear from Lemma \ref{crystal of product}.
            And then we know that  $ \mathcal{M}\left(Y_{p}(m) \right) \cdot\mathcal{M}\left(Y_{q}( 1 )\right)$
                  is a union of some connected components. Each connected component is isomorphic to some $ B \left( \lambda\right)(\lambda\in P_{+})$ since  the crystal $\mathcal{M}_{c}$ is normal
             by Corollary \ref{normality}.
           
           \PrfEnd

             \CorBgn \label{general result}
             For any $m_1,m_2, \mathellipsis m_l\in  \mathbb{ Z}$,  $p_1,p_2, \mathellipsis p_l
             \in I\, (l>0)$.
                  \begin{center}
             $ \mathcal{M}\left(Y_{p_1}( m_1) \right) \cdot  \mathcal{M}\left(Y_{p_2}(   
              m_2) \right) \cdot  \mathellipsis  \cdot  \mathcal{M}\left(Y_{p_l}( m_l) \right)$ ,
                \end{center}
               holds a crystal structure.
             \CorEnd

            \LemBgn  \label{decomposition}
            We have the following decomposition for the tensor product
            of fundamental crystals $B(\Lambda_{p})$ and $B(\Lambda_{q})$: 
              
           \begin{enumerate}
          \item For $p+q>n$,
           \begin{center}
            $B\left( \Lambda_{p} \right)\otimes  B\left( \Lambda_{q} \right)\cong  B\left( 
            \Lambda_{\max(p,q)}+\Lambda_{\min(p,q)} \right)\oplus  B\left( \Lambda_{{\max(p,q)}+1}
            +\Lambda _{\min(p,q)-1} \right)\oplus \mathellipsis
          \newline
          \oplus  B\left(\Lambda_{n}+\Lambda 
           _{p+q-n}\right)\oplus B\left(\Lambda_{p+q-n-1}\right).$
          \end{center}
          \item For $p+q\le n$,
            \begin{center}
           $ B\left( \Lambda_{p} \right)\otimes  B\left( \Lambda_{q} \right)\cong B\left( \Lambda_ 
           {\max(p,q)}+\Lambda_{\min(p,q)}  \right)\oplus  B\left( \Lambda_{{\max(p,q)}+1}+\Lambda 
           _{\min(p,q)-1} \right)\oplus \mathellipsis 
            \newline
           \oplus  B\left(\Lambda_{p+q-1}+\Lambda 
                   _{1}\right)\oplus  B\left( \Lambda_{p+q} \right).$
            \end{center}
                \end{enumerate}
           \LemEnd

          \PrfBgn
              cf. \cite{K}.
          \PrfEnd

         \LemBgn \label{use for Thm}
                    Any connected component in $ \mathcal{M}\left(Y_{p}( 
                    m) \right) \cdot\mathcal{M}\left(Y_{q}( 1 )\right)$
                    appears in $ \mathcal{M}\left(Y_{p}( m )
                    \right)\otimes  \mathcal{M}\left(Y_{q}( 1 )\right)$.
             \LemEnd

              \PrfBgn 
                    We may show that for each dominant weight $\lambda$ in $wt\left(( \mathcal{M}(Y_{p}\left( m\right))\otimes \mathcal{M}(Y_{q}\left( 1 \right))\right)$ we find $ B(\lambda)\subset  \mathcal{M}(Y_{p}\left( m \right))\,\otimes \, \mathcal{M}(Y_{q}\left( 1 \right))$,
                     since 
\begin{center}
                     $wt\left(( \mathcal{M}(Y_{p}\left( m\right))\cdot \mathcal{M}(Y_{q}\left( 1 \right))\right)
                     = wt \left( \mathcal{M}(Y_{p}\left( m \right))\otimes \mathcal{M}(Y_{q}\left( 1 \right))\right)$.
\end{center}
                     To find all the dominant weights in $ \mathcal{M}(Y_{p}\left( m \right))\otimes 
                      \mathcal{M}(Y_{q}\left( 1 \right))$, we may consider 
                     $B(\Lambda_{p})\otimes B(\Lambda_{q})$ as in  \cite{K}.
                     Let $(\,i_{1},i_{2},\mathellipsis,i_{p})$ be an element in $B(\Lambda_{p})$ 
                     like as a column tableau with the entries $i_{1},i_{2},\mathellipsis,i_{p}$ satisfying
                     $1\le i_{1}<i_{2}<i_{3}\mathellipsis<i_{p}\le n+1$.
                     Take $u
                     =( i_{1},\mathellipsis ,i_{p})\in B(\Lambda_{p})$
                     and $v=( j_{1},\mathellipsis ,j_{q})\in B(\Lambda_{q})$ whose weights are 
                    $wt( u)=\epsilon_{i_{1}}+\mathellipsis +\epsilon_{i_{p}}$,  
                     $wt( v )=\epsilon _{j_{1}}+\mathellipsis +\epsilon_{j_{q}}$,
                     where $\epsilon_{1}=\Lambda_{1},\, \epsilon_{i}=\Lambda_{i}-\Lambda_{i-1}\,(
                    i=1,2, \mathellipsis,n)$ and $\epsilon_{n+1}=-
                    \epsilon_{1}-\epsilon_{2}\mathellipsis-\epsilon_{n}$. Note that $\{\,\epsilon_{1}
                     \mathellipsis\epsilon_{n}\}\,$ are linearly independent.
                         Now  we assume $p\ge q$ without loss of generality.
                      Note that
                  \[
                  \Lambda_{k}=\epsilon_{1}+\mathellipsis +\epsilon_{k},\, \, \left( 
                   1\le k\le n\right).
                       \]
                           If  $wt\left( u \right)+ wt \left( v \right)$ is dominant, it should includes 
                            $\epsilon_{1}$  once or twice  and then we have the following cases:
         \vs{-.5ex}         \begin{enumerate}
                   \item \label{case 1}  If it includes $\epsilon_{1}$ once
                   then $wt\left( u \right)+wt\left( v \right)=\Lambda_{k}$  for some $k$ .
                   \item  \label{case 2} If it includes $\epsilon_{1}$ twice then $wt\left( u \right)+
                   wt\left( v \right)=\Lambda_{k}+\Lambda_{l}$ for some $k$ and $l$\,$(k\le l)$.
                           \end{enumerate}
                  Consider the case \ref{case 1}
                  \item First, we assume  $p+q>n$,
                   then
\[
                \Lambda_{k} = wt( u )+wt( v)=\epsilon_{i_{1}}+\mathellipsis 
                +\epsilon_{i_{p}}+\epsilon_{j_{1}}+\mathellipsis +\epsilon 
                  _{j_{q}}.
\]
                If neither $i_{p}$ nor $j_{q}$ is $n+1$, 
                some coefficients of $\epsilon_{l}\,(1\le l\le n)$ are 2, which means
                \begin{center} 
                $wt(u)+wt(v) \ne {\Lambda_{k}}$.
                 \end{center}
                If both $i_{p}$ and $j_{q}$ are $n+1$, it is trivial that 
                $wt\left(u\right)+wt\left(v\right)$ is never dominant.
                 \newline
                If  $i_{p}\ne n+1$ and $j_{q}=n+1$, then 
                $wt(u)$ and $wt(v)$ is
                  \begin{align}
                   \nonumber
                 wt( u )&=\epsilon_{i_{1}}+\mathellipsis +\epsilon_{i_{p}},\\
                       \nonumber
                  wt( v )&=\epsilon _{j_{1}}+\mathellipsis +\epsilon_{j_{q-1}}
                   -(\epsilon_{1}+\mathellipsis+\epsilon_{n}).
                    \end{align}
                  If $wt(u)+wt(v)=\Lambda_{k}=\epsilon_{1}+\mathellipsis+\epsilon_{k}$,
                  then we have
                \begin{align}
                 \epsilon_{1}+\mathellipsis+\epsilon_{k}=\epsilon_{i_{1}}+\mathellipsis 
                +\epsilon_{i_{p}}+\epsilon_{j_{1}}+\mathellipsis +\epsilon 
                  _{j_{q-1}}-(\epsilon_{1}+\mathellipsis 
                +\epsilon_{n}).
                 \end{align}
                  Then we obtain 
                   $k+n=p+q-1$ and then $k=p+q-n-1$. So, we get $q>k$. Indeed, if $q\le k=p+q-n-1$, then $0\le p-n-1$, which can not happen.
In this case $u=(i_{1},\mathellipsis, i_{p})=(1,\mathellipsis, p)$
                   and $v=(j_{1},\mathellipsis, j_{q-1}, j_{q})=(1,\mathellipsis ,k,p+1,\mathellipsis ,n,n+1)$ 
                  gives $wt(u)+wt(v)=\Lambda_{k}=\Lambda_{p+q-n-1}$. By Lemma 
                  \ref{decomposition}, the component $B(\Lambda_{p+q-n-1})$ appears
                  in $B(\Lambda_{p})\otimes B(\Lambda_{q})$.
                 \newline
                 The case $i_{p}=n+1$ and $j_{q}\ne n+1$ is done similarly.
                  Indeed, $(i_{1},\mathellipsis, i_{p-1})=(1,\mathellipsis, p-1)$
                  and $(j_{1},\mathellipsis, j_{q})=(1,\mathellipsis, k,p ,\mathellipsis, n)$
                  gives $wt(u)+wt(v)=\Lambda _{k}=\Lambda _{p+q-n-1}$.

                  Second, we assume $p+q\le n$. If one of $i_{p}$ or $j_{q}$ is $n+1$, 
                   $wt(u)+wt(v)$ is not dominant.
                   \newline
                   Thus, we know that neither $i_{p}$ nor $j_{q}$ is $n+1$.
                   The assumption $wt(u)+wt(v)=\Lambda_{k}$ means 
                    $p+q=k$.
                   \newline
                   Indeed, if $u=(i_{1},\mathellipsis, i_{P})=(1,\mathellipsis ,p)$
                   and $v=(j_{1},\mathellipsis, j_{q})=(p+1,\mathellipsis, p+q)$,
                  then we have $wt(u)+wt(v)=\Lambda_{p+q}$.
                    By Lemma  \ref{decomposition}, the component 
                   $B(\Lambda_{p+q})$ appears
                  in $B(\Lambda_{p})\otimes B(\Lambda_{q})$.

                 Now consider the case \ref{case 2}.
                 In this case, $\epsilon_{1}$ appears twice in 
                 $wt(u)+wt(v)$, which means neither $i_{P}$ nor $j_{q}$
                 is $n+1$. Then  we know that $p+q=k+l$.
                 \newline
                  First, assume $p+q>n$.
                 Now, we claim $k\le q\le p\le l$. Let us assume $k>q $. Then we have
                 \begin{align}  \label{eq6.2}
                 \Lambda_{k}+\Lambda_{l}=2\varepsilon_{1}+\mathellipsis +2\epsilon 
                  _{k}+\epsilon_{k+1}+\mathellipsis +\epsilon_{l}\, ,
                  \end{align}    \vs{-.5ex}
              and 
              \vs{-.5ex}   \begin{align} \label{eq6.3}
                    wt(u)+wt(v)=\epsilon_{i_{1}}+\mathellipsis 
                +\epsilon_{i_{p}}+\epsilon_{j_{1}}+\mathellipsis +\epsilon 
                  _{j_{q}},
                 \end{align}
                  which shows the number of coefficient 2 in \eqref{eq6.3}
                  is at most $q$ and then it contradicts with \eqref{eq6.2}.
                  Thus, we have $k\le q$ and $k\le q\le p\le l$
                  since $k+l=p+q$.
                   \newline
                  Such $(k,l)$ are:
                   $ (q,p),(q-1,p+1),(q-2,p+2),\mathellipsis ,(p+q-n,n)$.
                  All of them appear in the decomposition as in Lemma \ref{decomposition}.
                  Indeed, in this case $k=q-i$,\,$l=p+i$\,$(0\le i\le n-p)$,
                  $u=(i_{1},\mathellipsis, i_{p})=(1,\mathellipsis ,p)$
                   and $v=(j_{1},\mathellipsis, j_{q})=(1,\mathellipsis, q-i,p+1,\mathellipsis, p+i)$
                   gives $wt(u)+wt(v)=\Lambda_{p+i}+\Lambda_{q-i}$\,$(0\le i\le n-p)$.
  
                  Next, we assume $p+q\le n$. By arguing similarly to the previous case 
                  we find the possibility of $(k,l)$ are: 
                  $(q,p),(q-1,p+1),  \mathellipsis ,(1,p+q)$.
                   All of them appears in the decomposition as in Lemma \ref{decomposition}.
                   Indeed, in the case $k=q-i,\,\,l=p+i\,\,(0\le i\le q-1)$,
                     $u=(i_{1},\mathellipsis, i_{p})=(1,\mathellipsis, p)$
                   and $v=(j_{1},\mathellipsis ,j_{q})=(1,\mathellipsis ,q-i,p+1,\mathellipsis, p+i)$
                     gives $wt(u)+wt(v)=\Lambda_{p+i}+\Lambda_{q-i}$\,$(0\le i\le q-1)$.
\PrfEnd

   \vs{.25ex}
               \LemBgn \label{M=Yp(m)}
                     Let $M_{1}\in\mathcal{M}\left(Y_{p}( m )\right)$, 
                    $M_{2}\in\mathcal{M}\left(Y_{q}( 1) \right)$  
                     with $m\ge 1$ and $p,q\in \{1,\mathellipsis ,n\}\,$. 
                     If $M_{1}\cdot{M}_{2}$ is the highest weight vector then 
                     $M_{1}=Y_{p}\left( m \right)$.
                        \LemEnd

                 \PrfBgn
                The general form of $M_{1}$ is given in Corollary \ref{form of M},  
\[
                 M_{1}=X_{i_{1}}\left( p+m-1 \right)X_{i_{2}}\left( p+m-2 
                \right)\mathellipsis X_{i_{p}}(m)\,\,(1\le i_{1}<\mathellipsis <i_{P}\le n+1).
\]
                Assume that there is $j\in \{\,0,1,\mathellipsis,p-1\}\,$ such that
                $i_{j}+1<i_{j+1}$ (we set $i_{0}=0$) and take $j_{0}$ is the 
                biggest one among such $j'$s, which means there are gaps in the sequence
                $(i_{1}\mathellipsis i_{p})$ and 
                $j_{0}$ is the last gap. Note that $j_{0}<p$.
               Set $i_{j_{0}+1}=h$, 
                then there is 
\[
                X_{h}\left( p+m-{j_{0}}-1 \right)=\frac{Y_{h}(p+m-{j_{0}}-1)}{Y_{h-1}(p+m-{j_{0}})},
\]
               in $M_{1}$ and by the assumption that $M_{1}\cdot M_{2}$ is the highest, we should have the factor $Y_{h-1}(\alpha)$ in $ M_{2}$, and by the explicit form of $ M_{2}$, we can find the factors
\[
               \frac{Y_{h-1}(\alpha)}{Y_{h-2}(\alpha+1)}\cdot\frac{Y_{h'}(\alpha-1)}{Y_{h'-1}(\alpha)}\,\,
                 (h'>h),
\]
                  in $M_{2}$ and we obtain
             \begin{equation}
             \label{eq1}
             \alpha\ge p+m-j_{0}\,.
              \end{equation}
             We have $\,\frac{1}{Y_{h'-1}(\alpha)}\,$ in $M_{2}$, which means there is the factor 
              $Y_{h'-1}(\beta)$ in $M_{1}$ by the assumption that $M_{1}\cdot M_{2}$ is the highest.
              Since $h=i_{{j_{0}}+1}$ is the gap in $M_{1}$ and after $h$ the indices $i_{j}$'s are consecutive, 
               $h'-1$ should be $i_{p}$.

            Thus,  we have $Y_{h'-1}(\beta)=Y_{i_{p}}(m)$, and
         \vs{-1ex}    \begin{equation}
            \label{eq2}
            m\ge \alpha,   \vs{-.5ex}
            \end{equation}
            \vs{-1ex}      then from \eqref{eq1} and  \eqref{eq2}, we get
\[
             m\ge \alpha\ge p+m-{j_{0}}\ge m+1 ,
\]
            which is a contradiction.
            Then, we obtain $i_{j}=j$ (for all $j=1,\mathellipsis ,p$) and then $M_{1}=X_{1}(p+m-1)\mathellipsis
             X_{p}=Y_{p}(m)$.
\PrfEnd

\subsection{Decomposition of Monomial Product}

\LemBgn \label{dec}
                  Let $M$ be a monomial in $\mathcal {M}\left(Y_q(1)\right)$ and $p\in\{\,1,\mathellipsis,n\}\,$.
                 \begin{enumerate}
                \item If $wt(M)=\Lambda_{q}$, then
\vs{-2.75ex}
   \create{equation} \label{dec1}
\hs{-6em}M= X_{1}\left(q\right)\mathellipsis X_{q}\left(1\right)=Y_{q}(1).
\delete{equation}
                   \item If $max(1,p+q-n)\le i<\min(p, q)$ and $wt(M)=\Lambda_{p+q-i}+\Lambda_{i}-
                  \Lambda_{p}$, then
                \begin{align} \label{dec2}
               M=& X_{1}\left(q\right)\mathellipsis
              X_{i}\left(q-i+1\right)\cdot X_{p+1}\left(q-i\right)\mathellipsis X_{p+q-i}\left(1\right)\\
\nonumber
 =&Y_{i}(q-i+1)\cdot \frac{Y_{p+q-i}(1)}{Y_{p}(q-i+1)}.
                \end{align}
                \item If $p+q>n$ and $wt(M)=\Lambda_{p+q-n-1}-\Lambda_{p}$, then 
                \begin{align} \label{dec3}
                M=&X_{1}\left(q\right)\mathellipsis X_{p+q-n-1}\left(n-p+2\right)\cdot X_{p+1}\left(n-p+1\right)
               \mathellipsis X _{n+1}\left(1\right)\\
\nonumber
             =& \frac{Y_{p+q-n-1}(n-p+2)}{Y_{p}(n-p+2)}.
               \end{align}    
               \item If $p+q\le n$ and $wt(M)=\Lambda_{p+q}-\Lambda_{p}$, then 
              \vs{-4.75ex} \begin{align} \label{dec4}
             \hs{16em}   M=X_{p+1}\left(q\right)\mathellipsis X _{p+q}\left(1\right)
                  =\frac{Y_{p+q}(1)}{Y_{p}(q+1)}.
               \end{align}
               \end{enumerate}
\LemEnd

\PrfBgn
              It is known that all weight multiplicities of $B(\Lambda_{q})$ is free. And it is clear that the monomial
              \eqref{dec1},\eqref{dec2}, \eqref{dec3} and \eqref{dec4} are in  $\mathcal {M}\left(Y_q(1)\right)$  by Corollary \ref{form of M}, and have the corresponding weights.
\PrfEnd

\ThmBgn  \label{decomposition2}
                 For $M_{1}\in\mathcal{M}\left(Y_{p}( m )\right)$, 
                 $M_{2}\in\mathcal{M}\left(Y_{q}( 1) \right)$ with $m\ge 1$ and
                 $p,q\in \{1,\mathellipsis ,n\}$,  
                 assume that $M_{1}\cdot{M}_{2}$ is the highest weight vector if and only if the following
                 cases occur: 
                 \begin{enumerate}
                    \item  $wt(M)=\Lambda_{p}+\Lambda_{q}$, we have  $M_{1}=Y_{p}\left( m \right),\,
                    \,M_{2}=Y_{q}\left(1 \right)$.
                  \item $wt\left( M_{1}\cdot{M}_{2} \right)= \Lambda_{p+q-i}+\Lambda_{i}$
                  with\,\, $ \max {\left( 1,p+q-n \right)\le i\le \min(p,q)-1}$,
                  and \,$m\ge q-i+1$,  we have\,\,\,
                     $M_{1}=Y_{p}\left( m \right)$,

 $M_{2}=X_{1}\left( q 
                    \right)\mathellipsis X_{i}\left( q-i+1 \right)X_{p+1}\left( q-i 
                    \right)\mathellipsis X_{p+q-i}\left( 1 \right)=Y_{i}(q-i+1)\cdot \frac{Y_{p+q-i}(1)}{Y_{p}(q-i+1)}$.
                    \item  $wt\left( M_{1}\cdot{M}_{2} \right)=\Lambda_{p+q-n-1}$ with\,\,$ p+q>n$,
                    and \,$m\ge n-p+2$,  we have 
                    $M_{1}=Y_{p}\left( m \right)$, 
                    \newline
                    $M_{2}=X_{1}\left(q\right)\mathellipsis X_{p+q-n-1}\left( n-p+2\right) X_{p+1}
                    \left( n-p+1 \right)\mathellipsis X_{n+1}\left( 1 \right)=\frac{Y_{p+q-n-1}(n-p+2)}{Y_{p}(n-p+2)}$.
                    \item $wt\left( M_{1}\cdot{M}_{2} \right)=  \Lambda_{p+q}$ with \,\,$p+q\le n$, and 
                     \,$m\ge q+1$,  we have
                     $M_{1}=Y_{p}\left( m \right)$,
                      \newline
                     $M_{2}=X_{p+1}\left( q \right)\mathellipsis
                     \mathellipsis X_{p+q}\left( 1 \right)=\frac{Y_{p+1}(1)}{Y_{p}(q+1)}$.
                    \end{enumerate}
                    \ThmEnd

                   \PrfBgn
                    First, we assume that $M_{1}\cdot{M}_{2}$ is the highest weight vector then by 
                   Lemma \ref{M=Yp(m)}, $M_{1}=Y_{p}(m)$. We know by Lemma \ref{use for Thm} that 
                  any connected component $B(\lambda)$ in
                    $ \mathcal{M}\left(Y_{p}(m) \right)\cdot \mathcal{M}\left(Y_{q}( 1 )\right)$
                    appears in $ \mathcal{M}\left(Y_{p}( m )
                    \right)\otimes  \mathcal{M}\left(Y_{q}( 1 )\right)$. Suppose $p+q>n$. Thus, by Lemma \ref{decomposition} we can 
                    write
                     \begin{center}
                     $B(\lambda)\subset  B\left( 
                     \Lambda_{\max(p,q)}+\Lambda_{\min(p,q)} \right)\oplus  B
                     \left(\Lambda_{{\max(p,q)}+1}+\Lambda 
                      _{\min(p,q)-1} \right)\oplus \mathellipsis
                     \newline
                      \oplus  B\left(\Lambda_{n}+\Lambda 
                    _{p+q-n}\right)\oplus B\left(\Lambda_{p+q-n-1}\right).$
                       \end{center}
         \begin{enumerate}
                  \item
                 We know that if  $M_{1}=Y_{p}(m)$ and  $M_{2}=Y_{q}(1)$, then $M_{1}\cdot
                  M_{2}$ is always a highest weight vector in 
                   $ \mathcal{M}\left(Y_{p}(m) \right)\cdot \mathcal{M}\left(Y_{q}( 1 )\right)$
                  for any $m\ge 1$.
                \item 
                    In the previous case, we have already consider the case 
                   $ \Lambda_{\max(p,q)}+\Lambda_{\min(p,q)}$.
                    Then we may see the highest weight 
                     $\Lambda_{p+q-i}+\Lambda_{i}$  $( \max {\left( 1,p+q-n \right)\le i\le
                     \min(p,q)-1})$.
                       Then 
                       $wt(M_{2})= \Lambda_{p+q-i}+\Lambda_{i}-\Lambda_{p}$.
                      By Lemma \ref{dec},
                        the unique possibility of $M_{2}$ is 
                        \begin{center}   \vs{-.5ex}
                       $  M_{2}=Y_{i}(q-i+1)\cdot \frac{Y_{p+q-i}(1)}{Y_{p}(q-i+1)}$,   \vs{-1ex}
                      \end{center}
                          \vs{-.5ex}  then 
                       \vs{-.5ex}   \begin{align}
                       \nonumber
                       M_{1}\cdot M_{2}=Y_{p}\left(m\right)\cdot Y_{i}\left( q-i+1 \right)\cdot \frac{  Y_{p+q-i}\left( 1 
                         \right)}{Y_{p}\left(q-i+1\right)},
                      \end{align}
                        by the highest condition,
                      we get \,$m\ge q-i+1$.     
                    \item   If $wt\left( M_{1}\cdot{M}_{2} \right)=\Lambda_{p+q-n-1}$ with\,\,$ p+q>n
                    $, then 
                         the unique possibility of $M_{2}$ is 
                        \vs{-.5ex}  \begin{center}
                      $M_{2}=\frac{Y_{p+q-n-1}(n-p+2)}{Y_{p}(n-p+2)}$,   \vs{-.5ex}
                       \end{center}
                     \vs{-.5ex}      then
                       \begin{center}   \vs{-1ex}
                       $M_{1}\cdot M_{2}=Y_{p}\left(m\right)\cdot \frac{Y_{p+q-n-1}(n-p+2)}{Y_{p}(n-p+2)}$,
                        \end{center}
                         by the highest condition,
                         we get $m\ge n-p+2$.
                  \item    We assume $ p+q\le n$, then by Lemma \ref{use for Thm} the component $B(\lambda)$ in 
                         $ \mathcal{M}\left(Y_{p}(m) \right) \cdot\mathcal{M}\left(Y_{q}( 1 )\right)$
                         appears in $ \mathcal{M}\left(Y_{p}( m )
                          \right)\otimes  \mathcal{M}\left(Y_{q}( 1 )\right)$. Thus,   by Lemma \ref{decomposition} we can write
                         \begin{center}
                         $B(\lambda)\subset B\left( \Lambda_ 
                        {\max(p,q)}+\Lambda_{\min(p,q)}  \right)\oplus  B\left( \Lambda_{{\max(p,q)}+1}+\Lambda 
                        _{\min(p,q)-1} \right)\oplus \mathellipsis 
                        \newline
                        \oplus  B\left(\Lambda_{p+q-1}+\Lambda 
                        _{1}\right)\oplus  B\left( \Lambda_{p+q} \right).$
                        \end{center}
                         If $wt\left( M_{1}\cdot{M}_{2} \right)=\Lambda_{p+q}$, this case it done 
                        by similarly as the previous case.
 
                      On the contrary, if we assume  that we have the cases it is trivial that 
                      $M_{1}\cdot M_{2}$ is the highest.
  \end{enumerate}
                     \PrfEnd

                   \TxtBgn
                   For a statement $P$, let $\mathcal{L}(P)=\left\{\begin{array}{rcl}
                   1 & \mbox{if} & P\,\, is \,true,\\ 0 & \mbox{if} & P\,\, is \,false.
                   \end{array}\right.
                    $
                  \TxtEnd

                   \ThmBgn \label{smy.deco.}
                    \begin{enumerate}
                   \item If  $p+q>n$, we have 
                    \[
                   \mathcal{M}\left(Y_{p}(m) \right) \cdot\mathcal{M}\left(Y_{q}( 1 )\right)
                   \cong  B(\Lambda_{p}+\Lambda_{q})\,\oplus\bigoplus \limits_{i=\max(p
                    +q-n,q+1-m)}^{\min(p,q)-1}
                    B(\Lambda_{p+q-i}+\Lambda_{i})\oplus
                   \mathcal{L}\left(m\ge n-p+2\right)\cdot B(\Lambda_{p+q-n-1}).
                   \]
                
                   \item If $p+q\le n$, we have 
                   \[
                    \mathcal{M}\left(Y_{p}(m) \right)\cdot \mathcal{M}\left(Y_{q}(1)\right)
                   \cong  B(\Lambda_{p}+\Lambda_{q})\,\oplus\bigoplus \limits_{i=
                   \max(1,q+1-m)}^{\min(p,q)-1}
                   B(\Lambda_{p+q-i}+\Lambda_{i})\oplus
                   \mathcal{L}\left(m\ge q+1\right)\cdot B(\Lambda_{p+q}).
                  \]
                 \end{enumerate}
                
               \ThmEnd

               \PrfBgn 
                    By Theorem \ref{decomposition2}, we know that $ B(\Lambda_{p}+\Lambda_{q})$
                    always appears in $\mathcal{M}\left(Y_{p}(m) \right) \cdot\mathcal{M}\left(Y_{q}( 1 )
                   \right)$ for any $m\ge 1$.
                   For
                   $\,\max(1,p+q-n)\le i\le  \min(p,q)-1$,
                  $B(\Lambda_{p+q-i}+\Lambda{i})$
                 appears in $\mathcal{M}\left(Y_{p}(m) \right) \cdot\mathcal{M}\left(Y_{q}( 1 )
                   \right)$ if and only if $ i\ge q-m+1 .$

                   Then, we have the condition 
                   \[
                   \max(1,p+q-n,q+1-m)\le i \le {\min(p,q)-1}.
                   \]
                   If $p+q>n$, the component $ B(\Lambda_{p+q -n-1})$
                  appears in $\mathcal{M}\left(Y_{p}(m) \right) \cdot\mathcal{M}\left(Y_{q}( 1 )\right)$  if and
                   only if
                   $m\ge n-p+2$.
                    If $p+q\le n$, the component $ B(\Lambda_{p+q})$
                  appears in $\mathcal{M}\left(Y_{p}(m) \right) \cdot\mathcal{M}\left(Y_{q}( 1 )\right)$ if and  only if 
                $m\ge q+1$.
               \PrfEnd

\ExpBgn
                From Example  \ref{crystal graph}, we know the crystal graphs of 
                 $\mathcal{M}\left( Y_{2} (1)\right)$ and   $\mathcal{M}\left
                ( Y_{5} (m)\right)$ in $A_{5}$.
 \begin{enumerate}
                   \item 
               Let $p=5, q=2$, then $Y_{5}(m)\cdot Y_{2}(1)\in
               \mathcal{M}\left(Y_{5}(m) \right)\cdot \mathcal{M}\left(Y_{2}( 1 )\right)$,
               $p+q=7>n=5$, 
   \begin{enumerate}
                   \item if $m=1$, we have 
                    \begin{align}
                    \nonumber
                   \mathcal{M}\left(Y_{5}(1) \right)\cdot \mathcal{M}\left(Y_{2}( 1 )\right)
                  & \cong   B(\Lambda_{5}+\Lambda_{2})\,\oplus\bigoplus \limits_{i=
                     \max(2,2)}^{\min(5,2)-1}
                    B(\Lambda_{p+q-i}+\Lambda_{i})\oplus
                   \mathcal{L}\left(1\ge 2\right)\cdot B(\Lambda_{p+q-n-1})\\
                   \nonumber
                   & \cong B(\Lambda_{5}+\Lambda_{2})\,\oplus\bigoplus \limits_{i=
                     2}^{1} B(\Lambda_{5+2-i}+\Lambda_{i} ),\\
                   \nonumber
                   \mathcal{M}\left(Y_{5}(1) \right)\cdot \mathcal{M}\left(Y_{2}( 1 )\right)
                   &\cong   B(\Lambda_{5}+\Lambda_{2}),
                   \end{align}
                   \item if $m=2$, we have 
                   \begin{align}
                   \nonumber
                   \mathcal{M}\left(Y_{5}(2) \right)\cdot \mathcal{M}\left(Y_{2}( 1 )\right)
                  &\cong B(\Lambda_{5}+\Lambda_{2})\,\oplus\bigoplus \limits_{i=
                  \max(2,1)}^{\min(5,2)-1}
                    B(\Lambda_{p+q-i}+\Lambda_{i})\oplus
                   \mathcal{L}\left(2\ge 2\right)\cdot B(\Lambda_{p+q-n-1})\\
                   \nonumber
                    & \cong B(\Lambda_{5}+\Lambda_{2})\,\oplus\bigoplus \limits_{i=
                     2}^{1} B(\Lambda_{5+2-i}+\Lambda_{i} )\oplus B(\Lambda_{1}),\\
                   \nonumber
                   \mathcal{M}\left(Y_{5}(2) \right) \cdot\mathcal{M}\left(Y_{2}( 1 )\right)
                   &\cong   B(\Lambda_{5}+\Lambda_{2})\oplus B(\Lambda_{1}).
                   \end{align}
                  \end{enumerate}    
  \item           Let $p=q=2$, then $Y_{2}(m)\cdot Y_{2}(1)\in
                 \mathcal{M}\left(Y_{2}(m) \right) \cdot\mathcal{M}\left(Y_{2}( 1 )\right)$m
               $p+q=4\le n=5$, 
                   \begin{enumerate}
                   \item if $m=1$, we have 
                    \begin{align}
                    \nonumber
                     \mathcal{M}\left(Y_{2}(1) \right)\cdot \mathcal{M}\left(Y_{2}(1)\right)
                   &\cong B(\Lambda_{2}+\Lambda_{2})\,\oplus \bigoplus \limits_{i=
                   \max(1,2)}^{\min(2,2)-1}
                   B(\Lambda_{p+q-i}+\Lambda_{i})\oplus
                   \mathcal{L}\left(1\ge 3\right) \cdot B(\Lambda_{p+q}),\\
                    \nonumber
                    \mathcal{M}\left(Y_{2}(1) \right)\cdot \mathcal{M}\left(Y_{2}( 1 )\right)
                   &\cong B(2\Lambda_{2})\,\oplus \bigoplus \limits_{i=2}^{1}
                   B(\Lambda_{2+2-i}+\Lambda_{i}) \\
                    \nonumber
                   &\cong  B(2\Lambda_{2}),
                  \end{align}
                    \item if $m=3$, we have 
                    \begin{align}
                    \nonumber
                     \mathcal{M}\left(Y_{2}(3) \right)\cdot \mathcal{M}\left(Y_{2}(1)\right)
                   &\cong B(\Lambda_{2}+\Lambda_{2})\,\oplus \bigoplus \limits_{i=
                    \max(1,0)}^{\min(2,2)-1}
                   B(\Lambda_{p+q-i}+\Lambda_{i})\oplus
                   \mathcal{L}\left(3\ge 3\right)\cdot  B(\Lambda_{p+q}),\\
                      \nonumber
                    \mathcal{M}\left(Y_{2}(3) \right) \cdot\mathcal{M}\left(Y_{2}( 1 )\right)
                   &\cong B(2\Lambda_{2})\oplus \bigoplus \limits_{i=1}^{1}
                   B(\Lambda_{2+2-i}+\Lambda_{i}) \oplus  B(\Lambda_{4})\\
                      \nonumber
                   &\cong  B(2\Lambda_{2})\oplus B(\Lambda_{3}+\Lambda_{1}) 
                   \oplus  B(\Lambda_{4}).
                  \end{align}
                     \end{enumerate}
                     \item 
               Let $p=2, q=5$, then $Y_{5}(m)\cdot Y_{2}(1)\in
               \mathcal{M}\left(Y_{5}(m) \right)\cdot \mathcal{M}\left(Y_{2}( 1 )\right)$,
               $p+q=7>n=5$, 
   \begin{enumerate}
                   \item if $m=1$, we have 
                    \begin{align}
                    \nonumber
                   \mathcal{M}\left(Y_{2}(1) \right)\cdot \mathcal{M}\left(Y_{5}( 1 )\right)
                  & \cong   B(\Lambda_{2}+\Lambda_{5})\,\oplus\bigoplus \limits_{i=
                     \max(2,5)}^{\min(5,2)-1}
                    B(\Lambda_{p+q-i}+\Lambda_{i})\oplus
                   \mathcal{L}\left(1\ge 5\right)\cdot B(\Lambda_{p+q-n-1})\\
                   \nonumber
                   & \cong B(\Lambda_{2}+\Lambda_{5})\,\oplus\bigoplus \limits_{i=
                     5}^{1} B(\Lambda_{2+5-i}+\Lambda_{i} ),\\
                   \nonumber
                   \mathcal{M}\left(Y_{2}(1) \right)\cdot \mathcal{M}\left(Y_{5}( 1 )\right)
                   &\cong   B(\Lambda_{2}+\Lambda_{5}),
                   \end{align}
                   \item if $m=5$, we have 
                   \begin{align}
                   \nonumber
                   \mathcal{M}\left(Y_{2}(5) \right)\cdot \mathcal{M}\left(Y_{5}( 1 )\right)
                  &\cong B(\Lambda_{2}+\Lambda_{5})\,\oplus\bigoplus \limits_{i=
                  \max(2,1)}^{\min(5,2)-1}
                    B(\Lambda_{p+q-i}+\Lambda_{i})\oplus
                   \mathcal{L}\left(5\ge 5\right)\cdot B(\Lambda_{p+q-n-1})\\
                   \nonumber
                    & \cong B(\Lambda_{2}+\Lambda_{5})\,\oplus\bigoplus \limits_{i=
                     2}^{1} B(\Lambda_{2+5-i}+\Lambda_{i} )\oplus B(\Lambda_{1}),\\
                   \nonumber
                   \mathcal{M}\left(Y_{2}(5) \right) \cdot\mathcal{M}\left(Y_{5}( 1 )\right)
                   &\cong   B(\Lambda_{2}+\Lambda_{5})\oplus B(\Lambda_{1}).
                   \end{align}
  \end{enumerate}
      \end{enumerate}
\ExpEnd

\newpage

	\create{thebibliography}{9}
	\addcontentsline{toc}{section}{References}
		\vspace{-1ex}
           \bibitem{JS}\emph{ J.Hong.\,S-J.Kang, ``Introduction to Quantum Groups and Crystal Bases," American Mathematical Society (2002).}
            \bibitem {SJD}\emph{S-J.Kang,  J-A.Kim,  D-U.Shin, ``Monomial realization of crystal basas for special linear algebras," Journal of Algebra 
             247(2004)629-624.}
             \bibitem{MK}\emph{ M.Kashiwara, ``Realization of Crystals," in: Contemp. Math., vol. 325,Amer. Math. Soc., 2003, pp. 133-139.}
              \bibitem{K}\emph{T.Nakashima, ``Crystal base a generalization of the Littlewood-Richardson rule for the classsical Lie algebras,"  Commun. Math. Phys. 154 (1993), no. 2, 215-243.}
                \bibitem{N} \emph{H.\,Nakajima, ``$t$-Analogs of $q$-Characters of Quantum affine  Algebras of Type $A_{N}$, $D_{N}$,"  in: Contemp. Math., vol.325, Amer. Math. Soc., 2003,pp. 141-160.}
               \bibitem{MK-N}  \emph{M.Kashiwara, T.Nakashima, ``Crystal Graphs for Representations of the $q$-Analogue of Classical Lie Algebras,"  J. of Algebra, 165(1994), 295-345.}
                \bibitem{M}\emph{ M.Kashiwara, ``Crystal base and Littelmann's refined Demazure character formula,"  Duke Math. J. 71(1993)-839-858.}
                \bibitem{M1} \emph{M.Kashiwara, ``Crystalizing the $q$-analogue of universal enveloping 
                 algebras,"  Comm. Math. Phys. 133(1990), 249-260.}
                   \bibitem{M2}\emph{ M.Kashiwara, `` On crystal bases of the $q$-analogue of universal 
                enveloping  algebras,"  Duke Math. J. 63(1991), 456-516.}
	\delete{thebibliography}

\end{document}